\documentclass[letterpaper,10pt,oneside,onecolumn,reqno]{amsart}

\usepackage{amsmath, amssymb}
\usepackage{url}
\usepackage{enumerate}
\usepackage[all]{xy}
\usepackage[left=1in,top=1in,right=1in,bottom=1in]{geometry}
\usepackage{setspace}
\usepackage{hyperref}
\usepackage{wrapfig}
\usepackage{graphicx}
\usepackage{appendix}
\usepackage{color}
\usepackage[usenames,dvipsnames]{xcolor}

\onehalfspace
\newcommand{\C}{\mathcal C}

\newcommand{\EE}{\mathbb E}

\newcommand{\F}{\mathcal F}

\newcommand{\M}{\mathcal M}

\newcommand{\PP}{\mathbb P}

\newcommand{\Q}{\mathbb Q}
\newcommand{\R}{\mathbb R}

\newcommand{\Z}{\mathbb Z}

\renewcommand{\span}{\operatorname{span}}
\renewcommand{\div}{\operatorname{div}}

\newcommand{\SPD}{\operatorname{SPD}}
\newcommand{\Sym}{\operatorname{Sym}}
\newcommand{\Lip}{\operatorname{Lip}}

\newcommand{\deuc}{d_{\operatorname{Euc}}}
\newcommand{\dhaus}{d_{\operatorname{Haus}}}

\newcommand{\Rshape}{R_{\mathrm{shape}}}

\newcommand{\E}{\mathrm{e}}				
\newcommand{\D}{\mathrm{d}}				
\newcommand{\sD}{\, \mathrm{d}}				
\newcommand{\oo}{\infty}

\newcommand{\SO}{\operatorname{SO}}
\newcommand{\Var}{\operatorname{Var}}

\theoremstyle{definition}

\numberwithin{equation}{section}
\newtheorem{env_thm}{Theorem}[section]

\newtheorem{env_lem}[env_thm]{Lemma}

\newtheorem{env_pro}[env_thm]{Proposition}

\newtheorem{env_exa}[env_thm]{Example}

\newtheorem*{env_con}{Conjecture}

\title[Geodesics of Random Riemannian Metrics: Supplementary Material]{Geodesics of Random Riemannian Metrics: \\ Supplementary Material}

\author[T. LaGatta]{Tom LaGatta}
\address{Courant Institute of Mathematical Sciences \\ New York University \\ 251 Mercer St. \\ New York, New York 10012}
\email{tlagatta@gmail.com}

\author[J. Wehr]{Jan Wehr}
\address{Department of Mathematics\\ The University of Arizona \\ 617 N. Santa Rita Ave. \\ P.O. Box 210089 \\ Tucson, AZ 85721}
\email{wehr@math.arizona.edu}

\date{\today}

\begin{document}

\begin{abstract}
	This is supplementary material for the main \emph{Geodesics} article by the authors. In Appendix \ref{app_gaussian}, we present some general results on the construction of Gaussian random fields. In Appendix \ref{app_shape}, we restate our Shape Theorem from \cite{lagatta2009shape}, specialized to the setting of this article. In Appendix \ref{app_geomgeod}, we state some straightforward consequences on the geometry of geodesics for a random metric. In Appendix \ref{geombg}, we provide a rapid introduction to Riemannian geometry for the unfamiliar reader. In Appendix \ref{analytictools}, we present some analytic estimates which we use in the article. In Appendix \ref{proof_mo_lem}, we present the construction of the conditional mean operator for Gaussian measures. In Appendix \ref{fermiproof}, we describe Fermi normal coordinates, which we use in our construction of the bump metric.
\end{abstract}

\maketitle

\include{geodesics_main_mar13_content}

\setcounter{page}{2}

\appendix

\part{Supplemental Material (on arXiv)}
	
\section{Construction of Gaussian Random Fields} \label{app_gaussian}

	In this appendix, we construct Gaussian tensor fields on $\R^d$, which we use to generate random Riemannian metrics. A Gaussian $2$-tensor field $\xi_{ij}(x)$ is entirely defined by its mean $2$-tensor $m_{ij}(x) = \EE \xi_{ij}(x)$ and its covariance $4$-tensor $c_{ijkl}(x,y) = \EE \xi_{ij}(x) \xi_{kl}(y)$. Throughout we assume the mean tensor is zero, and that the covariance tensor is symmetric, stationary and isotropic, and compactly supported.  Such a Gaussian tensor field is the ``source of randomness'' for our random Riemannian metric, which we define pointwise by $g(x) = \varphi\big( \xi(x) \big)$, where $\varphi$ is a function which sends a symmetric $2$-tensor to a positive-definite one, and acts spectrally.
	
	It $c(x)$ is a symmetric, stationary and isotropic, and compactly supported covariance function, we may generate a covariance tensor by setting $c_{ijkl}(x,y) := c(|x-y|) \cdot \big( \delta_{ik} \delta_{jl} + \delta_{il} \delta_{jk} \big).$ It follows that $c_{iiii}(x,y) = 2 c(|x-y|)$ and $c_{ijij}(x,y) = c(|x-y|)$ if $i\ne j$, with all other components equal to zero.  It is not trivial that there exist Gaussian covariance functions $c$ satisfying these conditions.  We present a family of examples $c_d$ due to Gneiting \cite[Equation (17)]{gneiting2002csc} which are compactly supported and $6$-times differentiable at the origin.  
		
		\begin{env_exa}[Gneiting's covariance function] \label{exa_gneitingcov}
			Let $d \ge 1$, and choose any integers $\kappa \ge 3$ and $\nu \ge \tfrac{d+1}{2} + \kappa$.  Let $B$ denote the Beta function, and define
				\begin{equation} \label{cdk}
					c_d(r) := c_{\nu, \kappa}(r) := \frac{1}{B(2\kappa,\nu+1)} \int_r^1 u(u^2 - r^2)^{\kappa-1} (1-u)^{\nu} \sD u \end{equation}
			when $0 \le r \le 1$, and set $c_{\nu,\kappa}(r) = 0$ for $r \ge 1$.  Gneiting's function $c_d(r)$ is non-trivial, compactly supported, and $6$-times differentiable.\footnote{We remark that if $\nu < \tfrac{d+1}{2} + \kappa$, then Gneiting's function $c_{\nu,\kappa}(r)$ is not a Gaussian covariance function.}  
					
		For a single example that works in low dimensions ($2 \le d \le 9$), we set $\kappa = 3$ and $\nu = 8$.  In this case, formula \eqref{cdk} takes the explicit form $c_{d}(r) = -\frac{1}{5} (r-1)^{11} \left(5 + 55r + 239r^2 + 429 r^3\right)$.  The covariance functions for larger values of $d$ are of a similar polynomial character.
		\end{env_exa}
	
	Henceforth, let $\xi_{ij}(x)$ be a real-valued, stationary and isotropic Gaussian random field on $\R^d$ with mean zero and covariance tensor $c_{ijkl}(x,y) = \E \xi_{ij}(x) \xi_{kl}(y)$. Let $\Q$ be the law of the random field $\xi$ on $\Omega = C^2(\R^d, \Sym)$.

\section{The Shape Theorem} \label{app_shape}
	
	If $g$ is a random Riemannian metric, then $d_g$ is a random distance function, so $(\R^d, d_g)$ is a random metric space.  Let $B_g(t) = \{ x : d_g(0,x) \le t \}$ denote the ball of radius $t$ (with respect to this distance) centered at the origin in $\R^d$.  
	
	As a distance function, $d_g$ satisfies the triangle inequality on $\R^d$; 
 Using Liggett's version \cite{durrett1996probability} of Kingman's subadditive ergodic theorem \cite{kingman1968ets}, one easily sees that there exists some non-random constant $\mu \ge 0$ such that for each $v \in S^{d-1}$, we have that $\tfrac 1 r d_g(0, rv) \to \mu$ almost surely and in $L^1$.  Due to rotation-invariance of our model, the constant $\mu$ does not depend on the direction $v$.  
	
	A priori, there is no guarantee that $\mu > 0$ or that this convergence is uniform; however, both of these statements are ensured by the Shape Theorem \cite{lagatta2009shape}.  Heuristically, if $|x-y| \gg 1$, then $d_g(x,y) \sim \mu |x-y|$.  
	
	Let $B(r) = \{ x : |x| \le r \}$ denote the Euclidean ball of radius $r$ centered at the origin.  The Shape Theorem states that random Riemannian balls grow asymptotically like Euclidean balls:  for large $t$, $B_g(t) \sim B(t/\mu)$ almost surely.  We will formally define the measure $\PP$ in Section 2. The following theorem is an important consequence.
			
		\begin{env_thm}[Shape Theorem] \label{shapecor}
		Let $\PP$ be the measure on $(\Omega, \F)$ introduced in Section \ref{sect_rrm}.  This measure satisfies $\PP(\Omega_+) = 1$; let $g$ denote a random Riemannian metric with respect to $\PP$.  The following statements hold:
			\begin{enumerate}[a)]
				\item There exists a non-random constant $\mu > 0$ such that, with probability one, $\tfrac{1}{r} d_g(0,rv) \to \mu$ as $r \to \oo$, uniformly in the direction variable  $v \in S^{d-1}$.  The convergence also occurs in $L^2$, uniformly in $v$. 
				\item For all $\epsilon > 0$, with probability one, there exists a random constant $\Rshape = \Rshape(g)$ such that if $r \ge \Rshape$, then
					$$B(r - \epsilon r) \subseteq B_g(\mu r) \subseteq B(r + \epsilon r ).$$
				Equivalently, with probability one, the rescaled ball $\tfrac 1 t B_g(t)$ converges to the Euclidean ball $B(1/\mu)$ (in the Hausdorff topology on compact sets).  The Euclidean ball $B(1/\mu)$ is called the \emph{limiting shape} of the model.
				\item \label{shapecor_complete} With probability one, the Riemannian metric $g$ is geodesically complete.  Consequently, with probability one, for all $x$ and $y$ in $\R^d$, there is a finite, minimizing geodesic $\gamma$ connecting $x$ to $y$ such that $d_g(x,y) = L_g(\gamma)$.  With probability one, the topology on $\R^d$ generated by the metric $d_g(x,y)$ is complete.
			\end{enumerate}		
		\end{env_thm}
		\begin{proof}
			The constant $\mu$ is independent of the direction $v$ since the measure $\PP$ is rotationally-invariant.  Part (a) is Proposition 3.3 of \cite{lagatta2009shape}.  Part (b) is Theorem 3.1 of \cite{lagatta2009shape}.  Part (c) is Corollary 3.5 of \cite{lagatta2009shape}.  The equivalence of geodesic completeness and topological completeness is a consequence of the Hopf-Rinow theorem \cite{lee1997rmi}.
		\end{proof}
	
	In \cite{lagatta2009shape}, we called our model Riemannian first-passage percolation in analogy with lattice models of first-passage percolation.  In such models, one assigns a random \emph{passage-time distribution} $\tau$ to the bonds of the lattice $\Z^d$, and this generates a random metric $d_\tau$ on $\Z^d$ analogous to our random Riemannian metric $d_g$.  In fact, Kingman \cite{kingman1968ets} developed his subadditive ergodic theory to analyze the distance function of lattice FPP.  One can think of lattice first-passage percolation as a random perturbation of the flat Euclidean geometry of $\Z^d$.  For comprehensive recent surveys of FPP, we suggest the surveys by Howard \cite{howard2004mfp} and Blair-Stahn \cite{blair2007first}; the surveys by Kesten \cite{kesten1180arp, kesten1987pta} are older, but contain many technical details.
	
	To prove the Shape Theorem in our context, one tessellates $\R^d$ by unit cubes, and considers a dependent FPP model on the lattice formed by the centers of those cubes.  For each $z \in \Z^d$, we define $\Lambda_z = \sup |g|$, the maximum eigenvalue of the random Riemannian metric $g$ on the cube centered at $z$.  The random field $\Lambda_z$ induces a model of dependent FPP, and we use estimates from this setting to prove the Shape Theorem in the continuum.  We revisit these techniques in Section \ref{sect_proofoffrontiertimes}.
	
	In \cite{lagatta2009shape}, our proof was based on the robust energy-entropy method of mathematical physics, where one shows that an event occurs with extremely low probability over one particular lattice path (``high energy''), but sums this over all possible lattice paths starting at the origin (``high entropy'').  One proves convergence by showing that the ``energy'' beats the ``entropy'', which allows to apply the Borel-Cantelli lemma.  A large-deviation estimate is essential to this method: the number of lattice paths at the origin grows exponentially in $n$, the length of each path, so the probabilities must decay exponentially in $n$ for the arguments to work.  By adapting more carefully the proof of Cox and Durrett, one can weaken the assumptions to a finite moment bound.  For example, the assumption $\EE \min\{\Lambda_1, \cdots, \Lambda_{3^d}\}^{2\cdot 3^d+1} < \oo$ is sufficient for the Shape Theorem to hold. 
 
 	It is a classical result that geodesics can be dually interpreted as a Hamiltonian flow on the cotangent bundle $T^* \R^d$, with Hamilonian $H = g^{ij}(x) p_i p_j$. Armstrong and Souganidis \cite{armstrong2012stochastic,armstrong2011stochastic} have proved a general result on stochastic homogenization of Hamilton-Jacobi equations in stationary, ergodic environments, and our shape theorem is a special case of their theorem.

\section{The Geometry of Geodesics} \label{app_geomgeod}
	
	We continue with the assumption that $d=2$ in order to use Theorem \ref{POVthm} to say something about the plane geometry of geodesic curves.  Let $\beta$ be a probability measure on the tangent bundle $T\R^2$ which is absolutely continuous with respect to the Lebesque measure.  Let $(X,V)$ be chosen according to $\beta$, independently of the random metric $g$, and consider the geodesic $\gamma := \gamma_{X,V}$ with these random initial conditions.  The curve $\gamma$ is a random plane curve, and Proposition \ref{nonconstkappa} demonstrates one consequence of this randomness.  
	
		
	We say that a plane curve $\gamma$ \emph{contains a straight line segment} if the Euclidean acceleration $\tfrac{\ddot \gamma}{|\dot \gamma|^2}$ is constant on some interval.  More generally, we say that the curve $\gamma$ \emph{contains a circular arc} if the Euclidean normal acceleration $w(t) = \tfrac{1}{|\dot \gamma|^3} \langle \dot \gamma, \ddot \gamma \rangle$ is constant on some interval.

	\begin{env_pro} \label{nonconstkappa}
		Suppose that $d=2$.  With probability one, the geodesic $\gamma$ contains neither straight line segments nor circular arcs.
	\end{env_pro}
	\begin{proof}
		For any $s$, let $U_s \in \F$ denote the event that the turning $w$ is constant on an interval containing $s$:
			\begin{equation}
				U_s = \{ \mbox{$\exists ~\delta > 0$ s.t. $w$ is constant on the interval $(s-\delta,s+\delta)$} \}, \end{equation}
		and note that $U_s = \sigma_{-s} U_0$.  We must show that $\PP(\cup_s U_s) = 0$.  We first prove that $\PP(U_0) = 0$, then use a simple approximation argument and Theorem \ref{POVthm} to prove the lemma.
	
		The event $U_0$ implies that
			\begin{equation} \label{previous}
				0 = \dot w(0) = \tfrac{1}{|\dot\gamma|^6} \big( \langle \ddot \gamma, \ddot \gamma \rangle + \langle \dot \gamma, \dddot \gamma \rangle - \langle \dot \gamma, \ddot \gamma \rangle \cdot 3|\dot\gamma|^{2-1} \langle \dot\gamma, \ddot \gamma \rangle \big)
			\end{equation}
		which we can compute using the geodesic equation $\ddot \gamma^k = -\Gamma_{ij}^k \dot\gamma^i \dot \gamma^j$.  
		
		Since $\dot w(0)$ depends measurably on the random metric $g$, it is a random variable.  In particular, the elements of the terms of the right side of \eqref{previous} are rational functions of $|\dot\gamma(0)| = \dot\gamma^1(0) = 1/\sqrt{g_{11}(0)}$, the Christoffel symbols $\Gamma_{ij}^k(0)$ and their derivatives $\Gamma_{ij,l}^k(0)$.  These are all real-valued random variables with continuous distributions, hence the random variable $\dot w(0)$ also has a continuous distribution.  In particular, $\dot w(0) \ne 0$ almost surely.  Consequently,
			$$\PP(U_0) \le \PP\big( \dot w(0) = 0  \big) = 0.$$
			
		Now, let $s_n$ be a (non-random) countable dense sequence in $\R$.  If the event $U_s$ occurs for some (random) $s$, then it also must occur for some nearby $s_n$ (random $n$).  Thus
			$$\PP\left(\cup_{s} U_s \right) = \PP\left(\cup_{s_n} U_{s_n} \right) \le \sum_{n=0}^\oo \PP( \sigma_{s_n}^{-1} U_0) = 0,$$
		since the measures $\PP \circ \sigma_{s_n}^{-1}$ are all absolutely continuous with respect to the measure $\PP$ by Theorem \ref{POVthm}, and $\PP(U_0) = 0$.
	\end{proof}

\section{Riemannian Geometry Background} \label{geombg}

	We present a very rapid introduction to Riemannian geometry, working in coordinates on the fixed manifold $\R^d$.  For a more detailed introduction, see the book by Lee \cite{lee1997rmi}.
	
	Let $\SPD$ be the set of symmetric, positive-definite matrices, and let $\Omega_+ = C^2(\R^d, \SPD)$ be the space of $C^2$-smooth symmetric, positive-definite quadratic forms on $\R^d$.  Every element $g \in \Omega_+$ induces a \emph{Riemannian metric} on $\R^d$, that is, a smoothly varying inner product on the tangent bundle $T \R^d$.
	
	The flat Euclidean metric $\delta \in \Omega_+$ is the $2$-tensor field which is everywhere equal to the identity matrix:  $\delta_{ij}(x) = \delta_{ij}$, where the symbol on the right-side is the Kronecker $\delta$.  Let $\E_i$ be the standard basis vectors in $\R^d$.  We write each Riemannian metric $g \in \Omega_+$ in coordinates by $g_{ij}(x) = \langle \E_i, g(x) \E_j \rangle$, where the brackets denote the standard Euclidean inner product.
	
	If $v = v^i \E_i$ and $w = w^i \E_i$ are tangent vectors to $\R^d$ at some point $x$ (following the Einstein convention of summing over repeated upper and lower indices), the inner product of $w$ and $w$ with respect to $g \in \Omega_+$ is $\langle v, g(x) w \rangle = g_{ij}(x) v^i w^j$.
	
	Fix some $g \in \Omega_+$.  For a single tangent vector $v \in T_x \R^d$, we denote by $\|v\|_g = \sqrt{\langle v, g(x) v \rangle} = \sqrt{g_{ij}(x) v^i v^j}$ and $|v| = \sqrt{\langle v, v \rangle} = \sqrt{\delta_{ij} v^i v^j}$ the Riemannian and Euclidean lengths of $v$, respectively.  For a $C^1$-curve $\gamma : [a,b] \to \R^d$, we define the Riemannian arc length of $\gamma$ by $L_g(\gamma) = \int_a^b \| \dot \gamma(t) \|_g \sD t$.  We say that a curve is finite if it has finite Euclidean arc length; for our model, Theorem \ref{shapecor} implies that finite curves have finite Riemannian length (for almost every $g$).  The Riemannian distance between two points $x$ and $y$ is defined by $d_g(x,y) = \inf_\gamma L_g(\gamma),$ where the infimum is over all $C^1$-curves $\gamma$ connecting $x$ to $y$.

	We follow the Riemannian geometry convention and write the inverse of $g$ as $g^{ij}(x) = (g^{-1})_{ij}(x)$, so that $g^{ij} g_{jk}$ is the identity matrix $\delta^i_k$ (again, following the Einstein convention).  We denote derivatives of $g_{ij}$ by indices following a comma, so that $g_{ij,k} := \tfrac{\partial}{\partial x^k} g_{ij}$ and $g_{ij,kl} := \tfrac{\partial}{\partial x^l} \tfrac{\partial}{\partial x^k} g_{ij}$.

	The Riemannian metric $g$ induces a canonical covariant derivative, by way of the Levi-Civita connection.  The expression $\nabla_V W$ denotes the covariant derivative of a vector field $W$ along the vector field $V$.  Let $\E_i$ denote the standard basis vectors in $\R^d$; we use the same notation to denote the vector fields $\E_i(x) = \E_i$ for all $x \in \R^d$.  The Christoffel symbols are defined by writing the covariant derivative in coordinates: $\nabla_{\E_i} \E_j = \Gamma_{ij}^k \E_k$. 
	
	The Riemann curvature tensor $R = {R_{ijk}}^l(g,x)$ quantifies how much the covariant derivatives fail to commute:  $R(X,Y) Z = \nabla_X \nabla_Y Z - \nabla_Y \nabla_X Z$.  In this article, the only place we use the general Riemann curvature is Appendix \ref{fermiproof}, where we describe Fermi normal coordinates for general Riemannian manifolds.  We will mainly use the scalar curvature $K$.  In the two-dimensional case, the Riemann curvature tensor is determined by the scalar curvature (this is formula \eqref{2dimRiemanncurvature}).  
	
	We can express the Christoffel symbols $\Gamma_{ij}^k$ and the scalar curvature $K$ as polynomials in the metric and its first two derivatives:
		\begin{equation} \label{geoquantitiesdef}
			\Gamma_{ij}^k(g,x) = \tfrac 1 2 g^{km} \left(g_{im,j} + g_{mj,i} - g_{ij,k} \right) \qquad \mathrm{and} \qquad K(g,x) = g^{ij} \left( \Gamma_{ij,k}^k - \Gamma_{ik,j}^k + \Gamma_{ij}^k \Gamma_{kl}^l - \Gamma_{il}^k \Gamma_{kj}^l \right), \end{equation}
	where we evaluate the terms on the right side at the point $x$.

	A $C^2$-curve $\gamma$ is called a geodesic for the metric $g$ if its (covariant) acceleration is zero:  $\nabla_{\dot \gamma} \dot \gamma = 0$.  In coordinates, geodesics are solutions to the \emph{geodesic equation}
		\begin{equation} \label{geoeqn_geombg}
			\ddot \gamma^k(g,t) = -\Gamma_{ij}^k(g, \gamma(g,t)) \dot \gamma^i(g,t) \dot \gamma^j(g,t). \end{equation}
	We will often simplify our notation by writing $\ddot \gamma^k = -\Gamma_{ij}^k \dot\gamma^i \dot\gamma^j$.  The geodesic equation is the Euler-Lagrange equation for the Riemannian energy functional $L^2_g(\gamma) := \tfrac{1}{2} \int \|\dot\gamma\|^2_g$.

\section{Analytic Estimates} \label{analytictools}
	
	Throughout this article, it will be convenient to put a metric topology on the space of compact subsets of $\R^d$, in order to quantify what it means for two compact sets to be close to each other.    For any compact set $D \subseteq \R^d$, let $D^\epsilon = \{ x : \deuc(x, D) \le \epsilon \}$ denote the $\epsilon$-neighborhood of $D$, i.e. the set of all points which are Euclidean distance at most $\epsilon$ from $D$.  We define the Hausdorff metric on compact sets by
		\begin{equation}
			\dhaus(D, D') = \inf\{ \epsilon : \mbox{$D' \subseteq D^\epsilon$ and $D \subseteq (D')^\epsilon$} \}. \end{equation}
	That is, $\dhaus(D, D') \le \epsilon$ if and only if $D' \subseteq D^\epsilon$ and $D \subseteq (D')^\epsilon$.  Let $\C = \C(\R^d)$ be the space of compact sets equipped with the Hausdorff metric.  This is a complete metric space, and balls under this metric are compact.
	

	Next, we introduce the usual seminorms on differentiable tensor fields.  Let $\alpha = 0, 1$ or $2$.  For any compact $D \subseteq \R^d$, we define the usual $C^\alpha(D)$-seminorm of a symmetric quadratic form $\xi \in \Omega$ by
			$$\|\xi\|_{C(D)} = \sup_{x \in D} \max_{i,j} \{ |\xi_{ij}(x)| \}, \ \qquad \|\xi\|_{C^1(D)} = \sup_{x \in D} \max_{i,j,k} \{ |\xi_{ij}(x)|, |\xi_{ij,k}(x)| \}, \qquad \mathrm{etc.}$$ 
	We equip the space $\Omega$ with the topology generated by the $C(D)$-seminorms.  Consequently, the linear space $\Omega$ is a Fr\'echet space.  The space $\Omega_+$ of $C^2$-smooth Riemannian metrics is dense in the open cone $C(\R^d, \SPD) \subseteq \Omega$.

	For any function $f$ and compact $D \subseteq \R^d$, let $\Lip_D(f) = \limsup_{\epsilon \to 0} \sup_{D^\epsilon} |f(x) - f(y)|/|x-y|$ denote the Lipschitz constant of $f$ over an infinitesimal neighborhood of $D$.  We define the $C^{\alpha,1}(D)$-seminorms of a tensor field $\xi$ as the maximum of its $C^\alpha(D)$-seminorm and the $\Lip_D$-constant of the $\alpha$th derivatives of $\xi$:
			$$\|\xi\|_{C^{0,1}(D)} = \max_{i,j} \big\{ \|\xi\|_{C(D)}, \Lip_{D}(\xi_{ij}) \big\}, \qquad \|\xi\|_{C^{1,1}(D)} = \max_{i,j,k} \big\{ \|\xi\|_{C^1(D)}, \Lip_{D}(\xi_{ij,k}) \big\}, \qquad \mathrm{etc.}$$ 
	We use the same notation to denote the corresponding seminorms for scalar functions (and higher dimensional tensor fields).
		
	Let $\Omega^*$ denote the space of continuous linear functionals of $\Omega$.  Since each $f \in \Omega^*$ is a measurable function on the probability space $\Omega$, it is a random variable.  Of particular interest are the evaluation functionals $\delta^{kl}_x \in \Omega_+$, defined by $\delta^{ij}_x(\xi) = \xi_{ij}(x)$ for any $\xi \in \Omega$.  For any compact subset $D \subseteq \R^d$, let $\F_D^{(0)}$ denote the $\sigma$-algebra
		\begin{equation} \label{FD_def}
			\F_D^{(0)} = \bigcap_{\epsilon > 0} \sigma \big( \delta^{ij}_x : x \in D^\epsilon \big) \end{equation}
	generated by the information of a random tensor field on \emph{an infinitesimal neighborhood of} the set $D$.  Let $\F_D$ denote the completion of the $\sigma$-algebra $\F_D^{(0)}$.  
		
	Let $\|\cdot\|_D$ denote any of the $C^\alpha(D)$ or $C^{\alpha,1}(D)$-seminorms defined above for quadratic forms $\xi \in \Omega$.  Note that this definition also uses an infinitesimal neighborhood of $D$.  By definition, the seminorms $\|\cdot\|_D$ depend only on a tensor field  on \emph{an infinitesimal neighborhood} of the compact set $D$.  Consequently, 
		\begin{equation} \label{seminorms_FDmeasurable}
			\mbox{the seminorms $\|\cdot\|_D$ are $\F_D$-measurable.} \end{equation}			

	Recall that $Z_D(g) = \max\{ \|g\|_{C^{2,1}(D)}, \|g^{-1}\|_{C^{1,1}(D)} \}$ is the fluctuation functional for a random metric $g$ over the set $D$, defined in \eqref{Zdef}.
		
	\begin{env_lem} \label{ZD_cty}
		The map $(D, g) \mapsto Z_D(g)$ is jointly continuous in $D$ and $g$.  The random variable $g \mapsto Z_D(g)$ is $\F_D$-measurable.
	\end{env_lem}
	\begin{proof}
		The proof of the lemma follows easily from the statement,
			\begin{equation} \label{seminorms_jointlycts}
				\mbox{the map $(D, \xi) \mapsto \|\xi\|_{D}$ is jointly continuous.} \end{equation}
			
		To prove \eqref{seminorms_jointlycts}, we first fix a compact set $D_0 \in \C$ and a tensor field $\xi_0 \in \Omega$.  By construction, the map $D \mapsto \|\xi_0\|_{D}$ is continuous in the Hausdorff topology.  Let $\epsilon > 0$, and suppose that $D \subseteq D_0^\epsilon$, that $\big| \|\xi_0\|_D - \|\xi_0\|_{D_0} \big| \le \epsilon$, and that $\|\xi - \xi_0\|_{D_0^\epsilon} \le \epsilon$.  Then $\big| \|\xi\|_D - \|\xi_0\|_{D_0} \big| \le \|\xi - \xi_0\|_D + \big| \|\xi_0\|_D - \|\xi_0\|_{D_0} \big| \le \|\xi - \xi_0\|_{D_0^\epsilon} + \epsilon \le 2\epsilon$.
		
		The $\F_D$-measurability of $Z_D$ follows trivially from the $\F_D$-measurability of the seminorms $\|\cdot\|_D$, i.e., statement \eqref{seminorms_FDmeasurable}.
	\end{proof}
	
	The next lemma states that all the important geometric quantities are locally Lipschitz in the metric.  That is, the functions $\Gamma_{ij}^k(g,x)$,  $\Gamma_{ij,l}^k(g,x)$ and $K(g,x)$ are locally Lipschitz at each point $(g,x) \in \Omega_+ \times \R^d$.
			
	\begin{env_lem}[Uniform local Lipschitz estimate] \label{Lipest_lemma}
		Fix a compact set $D \subseteq \R^d$ and a metric $g \in \Omega_+$.  There exists constants $\epsilon$ and $L$ (both varying continuously in $D$ and $g$) such that if $\|g' - g\|_{C^{2,1}(D)} < \epsilon$, then 
			\begin{equation} \label{Lipest_Gamma}
				\big\| \Gamma(g', \cdot) - \Gamma(g, \cdot) \big\|_{C^{1,1}(D)} \le L \|g' - g\|_{C^{2,1}(D)} \end{equation}
		and
			\begin{equation} \label{Lipest_K}
				 \big\| K(g', \cdot) - K(g, \cdot) \big\|_{C^{0,1}(D)} \le L \|g' - g\|_{C^{2,1}(D)}. \end{equation}
	\end{env_lem}
	\begin{proof}
		Consider the space $\Omega_+ \times \R^d$ equipped with the product pseudometric $\|g' - g\|_{C^{2,1}(D)} \cdot |x' - x|$.  The terms $(g,x) \mapsto \big( g_{ij}(x), g_{ij,k}(x), g_{ij,kl}(x), g^{ij}(x), {g^{ij}}_{,k}(x) \big)$ are locally Lipschitz in $\Omega_+ \times \R^d$.  Since the Christoffel symbols, their derivatives and the scalar curvature $K$ are all polynomials in these terms, they are also locally Lipschitz.  Let $L$ be the largest such local Lipschitz constant, taken over all points in the compact set $D$.  Formulas \eqref{Lipest_Gamma} and \eqref{Lipest_K} follow.
	\end{proof}
			
	Recall that geodesics are solutions to the geodesic equation:  $\ddot \gamma^k = -\Gamma_{ij}^k \dot \gamma^i \dot \gamma^j$.  The geodesic equation is a second-order differential equation with locally-Lipschitz coefficients, so it is easily seen that geodesics exist, and are continuous in their initial conditions.  Let $\gamma_{x,v}$ be the unique geodesic with initial position $x$ and velocity $v$.
	
	\begin{env_lem} \label{geodesics_cts}
		The map $(x,v,g,t) \mapsto \gamma_{x,v}(g,t)$ is jointly continuous.
	\end{env_lem}	
	\begin{proof}
		Fix some metric $g \in \Omega_+$, and initial conditions $(x,v) \in (\R^d)^2$.  Let $T_*(g)$ be the blow-up time of $\gamma_g := \gamma_{x,v}(g,\cdot)$, and suppose that $t < T_*(g)$.  Let $D \subseteq \R^d$ be a compact set which contains the geodesic segment $\gamma_g|_{[0,t + \epsilon]}$ in its interior.
		
		The geodesic equation is a second-order ODE, with coefficients the Christoffel symbols.  Lemma \ref{Lipest_lemma} states that the Christoffel symbols $\Gamma_{ij}^k$ are locally Lipschitz functions of the metric and position.  Consequently, a theorem of smoothness of solutions of ODEs\footnote{A simple generalization of Theorem 31.8 of Arnol'd \cite{arnold1978ordinary}.} implies that
			\begin{equation} \label{Lipest_gamma}
				\sup_{s \in [0, t+\epsilon]} \big| \gamma_{x',v'}(g',s) - \gamma_{x,v}(g,s) \big| \le C \, \big\| \Gamma(g,\cdot) - \Gamma(g',\cdot) \big\|_{C^{0,1}(D)} \cdot |x - x'| \cdot |v - v'|, \end{equation}
		for some constant $C$ depending on $x$, $v$ and $g$.  By Lemma \ref{Lipest_lemma}, we have that $\| \Gamma(g,\cdot) - \Gamma(g', \cdot) \|_{C^{0,1}(D)} \le L \| g - g' \|_{C^{2,1}(D)}$, which implies the result.
	\end{proof}

\section{The Conditional Mean Operator $m_D$} \label{proof_mo_lem}

		For a compact set $D \subseteq \R^d$, let $D^\epsilon$ denote the $\epsilon$-neighborhood of $D$.  Recall that $\C$ is the space of compact sets of $\R^d$, equipped with the Hausdorff metric.  This means that that $\epsilon = \dhaus(D, D_0)$ is the smallest value of $\epsilon$ for which $D \subseteq D_0^\epsilon$ and $D_0^\epsilon \subseteq D$.
		
		If $A$ is any subset of $\Omega$, let $A^\perp \subseteq \Omega^*$ denote the annihilator of $A$:
			$$A^\perp = \{f \in \Omega^* : \mbox{$f(\xi) = 0$ for all $\xi \in A$} \}.$$

		\begin{env_lem} \label{kernel_eta}
			Fix a compact set $D \subseteq \R^d$.  On the subspace $K \Omega^*$ of $\Omega$, the restriction map $\eta_D$ has kernel $K(\eta_D^* X_D^*)^\perp$.
		\end{env_lem}
		\begin{proof}
			Let $f \in \Omega^*$.  For all $e \in \eta_D^* X_D^*$, the symmetry of $K$ implies that $e(\eta_D K f) = f(K \eta_D^* e)$.  Thus $f \in (\eta_D^* X_D^*)^\perp$ if and only if $Kf \in \ker \eta_D$.
		\end{proof}
		
		Lemma \ref{kernel_eta} implies that $\eta_D$ is injective on $K \eta^*_D X_D^*$, so the inverse map $\eta_D^{-1}$ is well-defined on $\eta_D K \eta^*_D X_D^*$.
		
		\begin{env_lem} \label{etaopnorm}
			For each compact $D \subseteq \R^d$, the linear map $\eta_D^{-1} : \eta_D K \eta^*_D X_D^* \to \Omega$ has operator norm $1$.
		\end{env_lem}
		\begin{proof}
		
			The evaluation functionals are dense in the space $X_D^*$, so the operator norm of $\eta_D^{-1}$ is given by
				\begin{equation} \label{opnormproof}
					\sup_{e \in X_D^*} \frac{\| K\eta_D^* e\|_{C(\R^d)}}{\|\eta_D^* e\|_{C(D)}} = \sup_{y \in D} \frac{\|c_y\|_{C(\R^d)}}{\|c_y\|_{C(D)}}. \end{equation}
			The covariance defines an inner product, hence it satisfies the Cauchy-Schwarz inequality.  For the stationary covariance function $c$, this implies that $c_y(x) = c(y,x) \le \sqrt{c(y,y) c(x,x)} = c_y(y)$.  Consequently, the numerator of \eqref{opnormproof} equals the denominator, and the ratio is constant and equal to $1$.
			
		\end{proof}
		

		\begin{proof}[Proof of Lemma \ref{mo_lem}] ~
		
		\emph{Proof of part \eqref{mo_identity}.}  By construction, $\eta_D \hat m_D = \eta_D \eta_D^{-1}$ is the identity operator on $\eta_D K \eta_D^* X_D^*$.  This space is dense in $X_D$ and the operator $\eta_D \hat m_D$ is continuous, so the identity property follows. \\
		
		\emph{Proof of part \eqref{mo_finrange}.}  By the finite-range dependence assumption on $c$, $c(y,x) = 0$ for any $y \in D$ and $|x-y| \ge 1$.  Since $K$ is the integral operator with kernel $c$, we have $K \eta_D^* e(x) = \int_D c(x,y) \sD \mu_e(y) = 0$ for any $x \notin D^1$ , so
			$$K \eta_D^* X_D^* \subseteq \{ \xi \in \Omega: \xi(x) = 0 \mathrm{~if~} x \notin D^1 \}.$$
		Clearly, $K \eta_D^* X_D^* = \hat m_D \big( \eta_D K \eta_D^* X_D^* \big)$.  The space $\eta_D K \eta_D^* X_D^*$ is dense in $X_D$ and $m_D$ is continuous, so this completes the proof. \\
		
		\emph{Proof of part \eqref{mo_measurability}.}  If $\xi(y) = \xi'(y)$ for all $y \in D$, then $\eta_D(\xi) = \eta_D(\xi')$ hence $\hat m_D \eta_D \xi = \hat m_D \eta_D \xi'$.  This proves that $\xi \mapsto \hat m_D \eta_D \xi$ is $\F_D$-measurable. \\
		
		\emph{Proof of part \eqref{mo_continuity}}  Fix a compact set $D_0 \subseteq \R^d$ and a function $\xi_0 \in \Omega$.  Write $\|\cdot\|$ for the seminorm $\|\cdot\|_{D^2}$.  Let $(D_n, \xi_n) \to (D, \xi)$, and choose $\epsilon \in (0,1)$ so that $D_n \subset D^\epsilon$ and $\|\xi_n - \xi\| \le \epsilon$.  In light of property \eqref{mo_finrange}, it suffices to prove that $\|m_{D_n}  \xi_n - m_D \xi\| \to 0$.
		
		Since the operator $K$ has dense image in $\Omega$ by Lemma \ref{K_lem}, we may choose $f \in KX^*$ so that
			\begin{equation}
				\| Kf - \xi \| < \epsilon, \end{equation}
		
		Using the triangle inequality, we calculate
			\begin{eqnarray}
				\|m_{D_n}  \xi_n - m_D  \xi\|  &\le&  \| m_{D_n}  \xi_n - m_{D_n} \xi\| + \|m_{D_n} \xi - m_{D_n} Kf\| \nonumber \\
				 && \qquad + \|m_{D_n} Kf - m_D  Kf\| + \|m_D  Kf - m_D  \xi \| \nonumber \\
				&\le& 3 \epsilon + \|m_{D_n}  Kf - m_D Kf\| \label{mo_continuity_proof2}
			\end{eqnarray}
		since the operators $m_{D_n}$ and $m_D$ all have operator norm $1$ by Lemma \ref{etaopnorm}.

		Using the identity \eqref{meta_iota}, we are able to transfer the problem into the setting of Hilbert spaces (as in Lemma \ref{monotoneKD}).  Combining the factorization $K = \iota \iota^*$, the identity \eqref{meta_iota}, and the fact that $\iota$ is a unitary map from $H$ to $X$, we have
			\begin{eqnarray}
				\|m_{D_n} Kf - m_D Kf\| &=& \|m_{D_n} \iota \iota^* f - m_D \iota \iota^* f\| = \|\iota \pi_{D_n} \iota^* f - \iota \pi_D \iota^* f \| \nonumber \\
				 &=& \|\pi_{D_n} \iota^* f - \pi_D \iota^* f \|_H, \label{mo_continuity_proof3}
			\end{eqnarray}
		where we denote the Hilbert-space norm by $\|\cdot \|_H$.
		
		We now introduce the idea of finite-dimensional projections.  For any $\delta > 0$ and any compact $D \subseteq \R^d$, define the $\delta$-mesh of $D$ by $\M_D^\delta := D \cap \delta \Z.$  Since the set $D$ is compact, the set $\M_D^\delta$ is finite.  Define the finite-dimensional subspaces
			$$H_D^\delta = \span\big\{ \iota^* \delta_x : x \in \M_D^\delta \big\} \subseteq H,$$
		and let $\pi_D^\delta$ denote the projection onto $H_D^\delta$ in $H$.  
		
		It is clear that $\overline{\cup H_D^\delta} = H_D$.  Consequently, for each fixed $f$, $\lim_{\delta\to 0} \pi_D^\delta \iota^* f = \pi_D \iota^* f$ in $H$.  Since $D_n \subseteq D^\epsilon$, we can choose $\delta$ uniformly:  there exists $\delta > 0$
			\begin{equation} \label{mo_continuity_proof4}
				\mbox{for all $n$, \quad $\| \pi_{D_n}^\delta \iota^* f - \pi_{D_n} \iota^* f\|_H < \epsilon$ \qquad and \qquad  $\| \pi_{D}^\delta \iota^* f - \pi_{D} \iota^* f\|_H < \epsilon$.} \end{equation}
				
		Since the sets $D_{n}$ converge to $D$ in the Hausdorff distance, for all sufficiently large $n$, the $\delta$-meshes $\M_{D_n}^\delta$ and $\M_D^\delta$ are equal.  Consequently, for all large $n$,
			\begin{equation} \label{mo_continuity_proof5}
				\pi_{D_n}^\delta = \pi_D^\delta. \end{equation}
				
		We are ready to complete the proof.  By applying the triangle inequality to \eqref{mo_continuity_proof3}, and using the estimates \eqref{mo_continuity_proof4} and \eqref{mo_continuity_proof5}, we have for all large $n$,
			\begin{eqnarray}
				\|\pi_{D_n} \iota^* f - \pi_D \iota^* f \|_H &\le& \|\pi_{D_n} \iota^* f - \pi_{D_n}^\delta \iota^* f\|_H + \|\pi_{D_n}^\delta \iota^* f - \pi_D^\delta \iota^* f \|_H + \|\pi_D^\delta \iota^* f - \pi_D \iota^* f \|_H \nonumber \\
				&\le& \epsilon + 0 + \epsilon = 2\epsilon. \label{mo_continuity_proof6}
			\end{eqnarray}
		By combining \eqref{mo_continuity_proof2}, \eqref{mo_continuity_proof3} and \eqref{mo_continuity_proof6}, we have proved that $\|m_{D_n} \xi_n - m_D \xi\| \le 5\epsilon$.  This completes the proof of Lemma \ref{mo_lem}.
		\end{proof}


	
	\section{Proof of Theorem \ref{fermilemma}, Existence of Fermi Normal Coordinates} \label{fermiproof}
				In Section 1.11 of \cite{poisson2004relativist}, Poisson derives the Fermi normal coordinates for the case of a pseudo-Riemannian metric in 4-dimensional spacetime.  The same analysis also works for Riemannian metrics in arbitrary dimension.  We focus on the general $d$-dimensional case here, then specialize to $d=2$ at the end of the proof to recover \eqref{gexpansion}.
				
				Let $\gamma(t)$ denote a geodesic along an arbitrary Riemannian manifold $(M,g)$.  Let $(\dot \gamma(t), n_2(t), \dots, n_d(t))$ be an orthonormal frame along $\gamma$.  Using the exponential map, define
					\begin{equation}
						\Phi_g(t,x^2, \dots, x^d) = \exp_{\gamma(t)}( x^i n_i(t) ). \end{equation}
				The coordinates $(t,x^2, \dots, x^d)$ are called \emph{Fermi normal coordinates}.  It is clear that in these coordinates, the geodesic is along the $t$-axis, and the Christoffel symbols vanish.  In the next lemma, we calculate the metric and its derivatives along the $t$-axis.  
				
				For notational convenience, we write symbols with more space, as with ${\Gamma^k}_{ij}$ instead of $\Gamma^k_{ij}$.  We also write subscripts with commas to denote partial derivatives, as with ${\Gamma^k}_{ij,l} := \tfrac{\partial}{\partial x^l} P{\Gamma^k}_{ij}.$  
				
				\begin{env_lem}
					\begin{equation} \label{gexpansion}
						\begin{array}{rcrcrcl}
							g_{11}(t,x) &=& 1 &-& R_{1k1l}(t) x^k x^l &+& O(x^3) \\
							g_{1j}(t,x) &=& &-& \tfrac 2 3 R_{1kjl}(t) x^k x^l &+& O(x^3) \\
							g_{ij}(t,x) &=& \delta_{ij} &-& \tfrac 1 3 R_{ikjl}(t) x^k x^l &+& O(x^3),
						\end{array}
					\end{equation}
				for $i$, $j$, $k$ and $l$ not equal to 1.
				\end{env_lem}
				\begin{proof}
							
				It follows easily from the definition of the Christoffel symbols that
					\begin{equation} \label{gijk}
						g_{ij,k} = g_{im} {\Gamma^m}_{kj} + g_{mj} {\Gamma^m}_{ik}. \end{equation}
				The vanishing of the Christoffel symbols on the geodesic $\gamma$ implies that $g_{ij,k} \equiv 0$ along $\gamma$.  To compute the second derivatives of $g_{ij}$, we will use the Riemann curvature tensor ${R^k}_{ijl}$, defined by
					\begin{equation} \label{riemanntensor}
						{R^k}_{ilj} = {\Gamma^k}_{ij,l} - {\Gamma^k}_{il,j} + {\Gamma^k}_{ml} {\Gamma^m}_{ij} - {\Gamma^k}_{mj} {\Gamma^m}_{il}, \end{equation}
				following the physics convention of ordering the indices.
			
				Since ${\Gamma^k}_{ij} \equiv 0$ along the geodesic,
					\begin{equation} \label{GammaFermi1}
						{\Gamma^k}_{ij,1} = 0, \end{equation}
				for any $i$, $j$ and $k$.  Plugging this into the definition \eqref{riemanntensor} of the Riemann curvature tensor gives
					\begin{equation} \label{GammaFermi2}
						{\Gamma^k}_{i1,l} = {R^k}_{il1}, \end{equation}
				for any $i$, $k$ and $l$.  The argument on page 23 of \cite{poisson2004relativist} implies that
					\begin{equation} \label{GammaFermi3}
						{\Gamma^k}_{ij,l} = -\tfrac 1 3 ({R^k}_{ijl} + {R^k}_{jil}), \end{equation}
				for any $k$, and for $i$, $j$ and $l$ not equal to 1.
				
				Since the metric is constant along $\gamma$, $g_{ij,1k} = 0$ for any $i$, $j$ and $k$.  Thus it suffices to calculate $g_{11,kl}$, $g_{1j,kl}$ and $g_{ij,kl}$ for $j$, $k$ and $l$ not equal to $1$. \newline
				
				Differentiating \eqref{gijk} and noting that the terms with Christoffel symbols vanish, we have
					\begin{equation} \label{gijklgeneral}
						g_{ij,kl} = g_{im} {\Gamma^m}_{kj,l} + g_{mj} {\Gamma^m}_{ik,l}, \end{equation}
				along the geodesic.  To calculate $g_{11,kl}$, we plug in the formula \eqref{GammaFermi2} for the first derivative of the Christoffel symbols to get
					\begin{equation} \label{g11kl}
						g_{11,kl} = 2 g_{1m} {\Gamma^m}_{k1,l} = 2g_{1m} {R^m}_{kl1} = 2R_{1kl1} = -2R_{1k1l}, \end{equation}
				where the last line follows from the symmetry $R_{1kl1} = R_{1k1l}$ of the Riemann tensor.  To calculate $g_{1j,kl}$, we apply both expressions \eqref{GammaFermi2} and \eqref{GammaFermi3} for the Christoffel symbols to \eqref{gijklgeneral} to get
					\begin{eqnarray}
						g_{1j,kl} &=& g_{1m} {\Gamma^m}_{kj,l} + g_{mj} {\Gamma^m}_{k1,l} = -\tfrac{1}{3} (R_{1kjl} + R_{1jkl}) + R_{jkl1} \nonumber \\
						&=& -\tfrac 1 3 R_{1kjl} + \tfrac 1 3 (R_{1ljk} + R_{1klj}) - R_{1ljk} \nonumber \\
						&=& -\tfrac 2 3 (R_{1kjl} + R_{1ljk}) \label{g1jkl}
					\end{eqnarray}
				where we use the symmetry $R_{jkl1} = - R_{1ljk}$, the Bianchi identity $R_{1jkl} = -R_{1ljk} - R_{1klj}$, and the symmetry $R_{1klj} = -R_{1kjl}$.
				
				By a similar argument, 
					\begin{equation} \label{gijklspec}
						g_{ij,kl} = -\tfrac 1 3 (R_{ikjl} + R_{ijkl} + R_{jikl} + R_{jkil}) = -\tfrac 1 3 (R_{ikjl} + R_{iljk}), \end{equation}
				where the middle two terms cancel by the symmetry $R_{ijkl} = -R_{jikl}$, and the last terms are equal by the symmetry $R_{jkil} = R_{iljk}$. \newline
				
				We now expand the metric $g(t,x)$ in a Taylor series around the point $(t,0)$, noting that $g_{ij}(t,0) = \delta_{ij}$, $g_{ij,k}(t,0) = 0$, and using the values \eqref{g11kl}, \eqref{g1jkl} and \eqref{gijklspec} for the second derivative $g_{ij,kl}(t,0)$ of the metric.  Formula \eqref{gexpansion} follows.
				\end{proof}

				In the case $d=2$, formula \eqref{gexpansion} takes a particularly simple form, since the Riemann curvature tensor is determined by the scalar curvature $K(t)$ via the following identity:
					\begin{equation} \label{2dimRiemanncurvature}
						R_{1212}(t) = \tfrac 1 2 K(t) \det g = \tfrac 1 2 K(t) (g_{11} g_{22} - g_{12}^2). \end{equation}
				Applying this, we have $R_{1212} = \tfrac 1 2 K(t)$, and the terms with $R_{1222}$ and $R_{2222}$ vanish by the symmetries of the curvature tensor, so
					\begin{equation} \label{gexpansion2}
						g_{11}(t,x) = 1 - \tfrac 1 2 K(t) x^2 + O(x^3), \qquad g_{12}(t,x) = O(x^3), \qquad \mathrm{and} \qquad g_{22}(t,x) = 1 + O(x^3). \end{equation}




\end{document}